\documentclass[10pt]{article}
\usepackage{latexsym}
\usepackage{amsfonts}
\usepackage{enumerate}
\usepackage{multicol}

\title{On the maximum size of a $(k,l)$-sum-free subset of an abelian group \\[.4in]}

\author{B\'{e}la Bajnok \\[.1in] Department of Mathematics, Gettysburg College \\
Gettysburg, PA 17325-1486 USA \\E-mail:  bbajnok@gettysburg.edu \\[.4in]}

\date{September 25, 2007}

\newtheorem{thm}{Theorem}

\newtheorem{lem}[thm]{Lemma}

\newtheorem{prop}[thm]{Proposition}

\newtheorem{que}{Question}

\begin{document}

\maketitle

\begin{abstract}

A subset $A$ of a given finite abelian group $G$ is called $(k,l)$-sum-free if the sum of $k$ (not necessarily distinct) elements of $A$ does not equal the sum of $l$ (not necessarily distinct) elements of $A$.  We are interested in finding the maximum size $\lambda_{k,l}(G)$ of a $(k,l)$-sum-free subset in $G$.  

A $(2,1)$-sum-free set is simply called a sum-free set.  The maximum size of a sum-free set in the cyclic group $\mathbb{Z}_n$ was found almost forty years ago by Diamanda and Yap; the general case for arbitrary finite abelian groups was recently settled by Green and Ruzsa.  Here we find the value of $\lambda_{3,1}(\mathbb{Z}_n)$.  More generally, a recent paper of Hamidoune and Plagne examines $(k,l)$-sum-free sets in $G$ when $k-l$ and the order of $G$ are relatively prime; we extend their results to see what happens without this assumption.        
\end{abstract}

\noindent 2000 Mathematics Subject Classification:  \\ Primary: 11P70; \\ Secondary: 05D99, 11B25, 11B75, 20K01.

\noindent Key words and phrases: \\ Sum-free sets, $(k,l)$-sum-free sets, Kneser's Theorem, arithmetic progressions.

\thispagestyle{empty}
\pagestyle{myheadings}
\markright{$(k,l)$-sum-free sets}

\section{Introduction}

Throughout this paper, we let $G$ be a finite abelian group of order $n>1$, written in additive notation; $v$ will denote the exponent (i.e. largest order of any element) of $G$. 

For subsets $A$ and $B$ of $G$, we use the standard notations $A+B$ and $A-B$ to denote the set of all two-term sums and differences, respectively, with one term chosen from $A$ and one from $B$.  If, say, $A$ consists of a single element $a$, then we simply write $a+B$ and $a-B$ instead of $A+B$ and $A-B$.  For a positive integer $h$ and a subset $A$ of $G$, the set of all $h$-term sums with (not necessarily distinct) elements from $A$ will be denoted by $hA$.   

Let $k$ and $l$ be distinct positive integers.  A subset $A$ of $G$ is called a \emph{$(k,l)$-sum-free set} in $G$ if $$kA \cap lA = \emptyset;$$ or, equivalently, if $$0 \not \in kA - lA.$$ Clearly, we may assume that $k>l$.  We are interested in determining the maximum possible size $\lambda_{k,l}(G)$ of a $(k,l)$-sum-free set in $G$.

A $(2,1)$-sum-free set is simply called a sum-free set.  The value of $\lambda_{2,1}(\mathbb{Z}_n)$ was determined by Diamanda and Yap \cite{DiaYap:1969a} in 1969.  It can be proved (see also \cite{WalStrWal:1972a}) that

\begin{equation} \label{(2,1)bounds}
\max_{ d | v} \left\{  \left \lfloor \frac{d+1}{3} \right \rfloor \cdot \frac{n}{d}  \right \}  \leq \lambda_{2,1}(G) \leq \max_{ d | n} \left\{  \left \lfloor \frac{d+1}{3} \right \rfloor \cdot \frac{n}{d}  \right \},
\end{equation}
which for cyclic groups immediately implies the following.  

\begin{thm} [Diamanda and Yap \cite{DiaYap:1969a}] \label{(2,1)cyclic}

The maximum size $\lambda_{2,1}(\mathbb{Z}_n)$ of a sum-free set in the cyclic group of order $n$ is given by
$$\lambda_{2,1}(\mathbb{Z}_n)=\max_{ d | n} \left\{  \left \lfloor \frac{d+1}{3} \right \rfloor \cdot \frac{n}{d}  \right \} = \left\{ 
\begin{array}{cl}
\frac{p+1}{p} \cdot \frac{n}{3} & \mbox{if    } n \mbox{   is divisible by a prime    } p \equiv 2  \mbox{    \rm{(mod}     3)   } \\
& \mbox{   and } p  \mbox{   is the smallest such prime;} \\ \\
\left \lfloor \frac{n}{3} \right \rfloor & \mbox{otherwise.    }  \\
\end{array}
\right.$$
\end{thm}     

The problem of finding $\lambda_{2,1}(G)$ for arbitrary $G$ stood open for over 35 years.  In a recent breakthrough paper, Green and Ruzsa \cite{GreRuz:2005a} proved that, as it has been conjectured, the value of $\lambda_{2,1}(G)$ agrees with the lower bound in (\ref{(2,1)bounds}):

\begin{thm} [Green and Ruzsa \cite{GreRuz:2005a}] \label{(2,1)all}
The maximum size $\lambda_{2,1}(G)$ of a sum-free set in $G$ is
$$\lambda_{2,1}(G)= \lambda_{2,1}(\mathbb{Z}_{v}) \cdot  \frac{n}{v} = \max_{ d | v} \left\{  \left \lfloor \frac{d+1}{3} \right \rfloor \cdot \frac{n}{d}  \right \}.$$
\end{thm}

As a consequence, we see that $$\frac{2}{7}n \leq \lambda_{2,1}(G) \leq \frac{1}{2}n$$ for every $G$, with equality holding in the lower bound when $v =7$ and in the upper bound when $v$ (iff $n$) is even. 

Now let us consider other values of $k$ and $l$.  In Section 2 of this paper we generalize (\ref{(2,1)bounds}), and prove the following.

\begin{thm} \label{boundgeneral}
The maximum size $\lambda_{k,l}(G)$ of a $(k,l)$-sum-free set in $G$ satisfies   
$$\max_{d | v} \left\{ \left( \left \lfloor \frac{d-1 - \delta(d)}{k+l} \right \rfloor +1  \right) \cdot \frac{n}{d}  \right\} \leq \lambda_{k,l}(G) \leq \max_{d | n} \left\{ \left( \left \lfloor \frac{d-2}{k+l} \right \rfloor +1  \right) \cdot \frac{n}{d} \right\},$$
where $\delta(d) = \mathrm{gcd} (d, k-l)$.

\end{thm}

Note that for $(k,l)=(2,1)$ Theorem \ref{boundgeneral} yields (\ref{(2,1)bounds}).  Note also that, if $k-l$ is not divisible by $v$, then $\delta(v) = \mathrm{gcd} (v, k-l) \leq v/2$; in particular, $$\lambda_{k,l}(G) \geq \frac{n}{2(k+l)} > 0.$$  If, on the other hand, $k-l$ is divisible by $v$, then clearly $\lambda_{k,l}(G)=0$, since for any $a \in G$ we have $ka = la$.

Let us now consider cyclic groups.  When $G \cong \mathbb{Z}_n$ and $n$ and $k-l$ are relatively prime, then Theorem \ref{boundgeneral} gives 
\begin{equation} \label{relprime}
\lambda_{k,l}(\mathbb{Z}_n) = \max_{d | n} \left\{ \left( \left \lfloor \frac{d-2}{k+l} \right \rfloor +1  \right) \cdot \frac{n}{d} \right\}.
\end{equation}
This result was already established by Hamidoune and Plagne in \cite{HamPla:2003a}.  Their method was based on a generalization of Vosper's Theorem \cite{Vos:1956a} on critical pairs where arithmetic progressions, that is, sets of the form $$A=\{ a, a+d, \dots, a+ c \cdot d \}$$ play a crucial role.  In particular, Hamidoune and Plagne proved that, if $G \cong \mathbb{Z}_n$ and $n$ and $k-l$ are relatively prime, then   
 \begin{equation} \label{alpha1}
\lambda_{k,l}(\mathbb{Z}_n) = \max_{d | n} \left\{  \alpha_{k,l}(\mathbb{Z}_d) \cdot \frac{n}{d} \right\},
\end{equation}
 where $\alpha_{k,l}(\mathbb{Z}_n)$ is the maximum size of a $(k,l)$-sum-free arithmetic progression in $\mathbb{Z}_n$.  Hamidoune and Plagne deal only with the case when $n$ and $k-l$ are relatively prime; as they point out, ``in the absence of this assumption, degenerate behaviors may appear'', and we concur with this assessment.  Nevertheless, we attempt to treat the general case; in Section 3 of this paper we prove that (\ref{alpha1}) remains valid even without the assumption that $n$ and $k-l$ are relatively prime:

\begin{thm} \label{cyclic}

For arbitrary positive integers $k$, $l$, and $n$ we have

$$\lambda_{k,l}(\mathbb{Z}_n) = \max_{d| n} \left\{ \alpha_{k,l}(\mathbb{Z}_d) \cdot \frac{n}{d} \right\}.$$

\end{thm}

Let us now move on to general abelian groups.  Hamidoune and Plagne conjecture in \cite{HamPla:2003a} that $$\lambda_{k,l}(G)= \lambda_{k,l}(\mathbb{Z}_{v}) \cdot  \frac{n}{v}$$ holds when $n$ and $k-l$ are relatively prime.  They prove this assertion with the additional assumption that at least one prime divisor of $v$ is not congruent to 1 (mod $k+l$).  We generalize this result for the case when $n$ and $k-l$ are not necessarily relatively prime:  

\begin{thm} \label{(k,l)general}
As before, for a positive integer $d$, we set $\delta(d)=\mathrm{gcd}(d,k-l)$.  If $v$ possesses at least one divisor $d$ which is not congruent to any integer between 1 and $\delta(d)$ (inclusive) (mod $k+l$), then $$\lambda_{k,l}(G)= \lambda_{k,l}(\mathbb{Z}_{v}) \cdot  \frac{n}{v}.$$ 
\end{thm}

We closely follow some of the fundamental work of Hamidoune and Plagne in \cite{HamPla:2003a}; in fact, Section 3 of this paper can be considered an extention of \cite{HamPla:2003a} for the case when $n$ and $k-l$ are not assumed to be relatively prime.

In Section 4 we employ Theorem \ref{cyclic} to establish the value of $\lambda_{3,1}(\mathbb{Z}_n)$ explicitly.  As an analogue to Theorem \ref{(2,1)cyclic} we prove the following.

\begin{thm}  \label{(3,1)cyclic}
The maximum size $\lambda_{3,1}(\mathbb{Z}_n)$ of a $(3,1)$-sum-free set in the cyclic group of order $n$ is given by
 $$\lambda_{3,1}(\mathbb{Z}_n)=\max_{ \begin{array}{c} d | n \\ d \not \equiv 2 \mbox{   } (\mathrm{  mod  } \mbox{   }  4) \end{array} } \left\{  \left \lfloor \frac{d+2}{4} \right \rfloor \cdot \frac{n}{d}  \right \} = \left\{ 
\begin{array}{cl}
\frac{p+1}{p} \cdot \frac{n}{4} & \mbox{if    } n \mbox{   is divisible by a prime    } p \equiv 3  \mbox{    \rm{(mod}     4)   } \\
& \mbox{   and } p  \mbox{   is the smallest such prime;} \\ \\
\left \lfloor \frac{n}{4} \right \rfloor & \mbox{otherwise.    }  \\
\end{array}
\right.$$
\end{thm}

As a consequence, we see that $$\frac{1}{5}n \leq \lambda_{3,1}(\mathbb{Z}_n) \leq \frac{1}{3}n,$$ with equality holding in the lower bound when $n \in \{5,10\}$ and in the upper bound when $n$ is divisible by 3.   

In our final section, Section 5, we provide some further comments and discuss several open questions about $(k,l)$-sum-free sets.

\section{Bounds for the size of maximum $(k,l)$-sum-free sets}

In this section we prove Theorem \ref{boundgeneral}.

We will use the following easy lemma.

\begin{lem} \label{stabil}

Suppose that $A$ is a maximal $(k,l)$-sum-free set in $G$.  Let $K$ denote the stabilizer subgroup of $kA$.  Then 

(i) $k(A+K)=kA$;

(ii) $A+K$ is a $(k,l)$-sum-free set in $G$;

(iii) $A+K=A$;

(iv) $A$ is the union of cosets of $K$.

\end{lem}

\emph{Proof.}  (i) The inclusion $kA \subseteq k(A+K)$ is obvious.  Suppose that $a_1, \dots, a_k \in A$ and $h_1, \dots, h_k \in K$.  Then $$(a_1+ \cdots +a_k) + (h_1+ \cdots +h_k) \in kA,$$ so $k(A+K) \subseteq kA$.

(ii) Suppose, indirectly, that $$k(A+K) \cap l(A+K) \not = \emptyset;$$ by (i) this implies $$kA \cap l(A+K) \not = \emptyset.$$  Then we can find elements $a_1, \dots , a_k \in A$, $a'_1, \dots , a'_l \in A$, and $h_1, \dots , h_l \in K$ for which $$a_1+ \cdots + a_k = a'_1+ \cdots +a'_l + h_1+ \cdots +h_l.$$  But $$ a'_1+ \cdots +a'_l = a_1+ \cdots + a_k -h_1- \cdots -h_l \in kA,$$ and this contradicts the fact that $A$ is $(k,l)$-sum-free.  

(iii) Since $A \subseteq A+K$ and $A$ is a maximal $(k,l)$-sum-free set in $G$, by (ii) we have $A+K=A$.

(iv) We need to show that for any $a \in A$, we have $a+K \subseteq A$.  But $a+K \subseteq A+K$, so the claim follows from (iii).  $\quad \Box$

For the upper bound in Theorem \ref{boundgeneral}, we need the following result which is essentially due to Kneser.

\begin{thm} [Kneser \cite{Kne:1956a}; see Theorem 4.4 in \cite{Nat:1996a}] \label{Kneser}

Suppose that $A$ is a non-empty subset of $G$ and, for a given positive integer $h$, let $H$ be the stabilizer of $hA$.  Then we have $$|hA| \geq h \cdot |A|- (h-1) \cdot |H|.$$  

\end{thm} 

\emph{Proof of the upper bound in Theorem \ref{boundgeneral}.}  Let $A$ be a $(k,l)$-sum-free set in $G$ with $|A|=\lambda$; then we have
$$kA \cap lA = \emptyset$$ and therefore 
\begin{eqnarray} \label{emptyset}
n \geq |kA| + |lA|.
\end{eqnarray}  
As before, let $K$ and $L$ be the stabilizer subgroups of $kA$ and $lA$, respectively.  Then, by Theorem \ref{Kneser}, we have 
$$|kA| \geq k \cdot |A|- (k-1) \cdot |K|$$ and $$|lA| \geq l \cdot |A|- (l-1) \cdot |L|;$$ thus, from (\ref{emptyset}) we get $$n \geq (k+l) \cdot |A| - (k-1) \cdot |K| - (l-1) \cdot |L|.$$
Without loss of generality we can assume that $|K| \geq |L|$, so
$$n \geq (k+l) \cdot |A| - (k+l-2) \cdot |K|$$ or 
$$\frac{|A|}{|K|} \leq \frac{1}{k+l} \cdot \left(\frac{n}{|K|} + (k+l-2) \right).$$
Now $|A|=\lambda$; in particular, $A$ is maximal, so by Lemma \ref{stabil} (iv), $\frac{|A|}{|K|}$ must be an integer.  Therefore, with $d$ denoting the index of $K$ in $G$, we get $$\frac{\lambda}{n/d} \leq \left \lfloor \frac{1}{k+l} \cdot \left(d + k+l-2 \right) \right \rfloor,$$ from which our claim follows.       
$\quad \Box$

\begin{prop} \label{delta}
Let $d$ be a positive integer, and set $\delta(d) = \mathrm{gcd} (d, k-l)$.  Suppose that $c$ is a positive integer for which $$(k+l) \cdot c \leq d-1- \delta(d).$$  Then there exists an element $a \in \mathbb{Z}_d$ for which the set $$A =\{a,a+1,a+2, \dots, a+c\}$$ is a $(k,l)$-sum-free in $\mathbb{Z}_d$ of size $c+1$.
\end{prop}

\emph{Proof.} By the Euclidean Algorithm, we have unique integers $q$ and $r$ for which $$l \cdot c = \delta(d) \cdot q -r$$ and $1 \leq r \leq \delta(d)$.  We also know the existence of integers $u$ and $v$ for which $$\delta(d) = (k-l) \cdot u + d \cdot v.$$  Now set $a=u \cdot q$.  We will show that $$A =\{a,a+1,a+2, \dots, a+c\}$$ is a $(k,l)$-sum-free in $\mathbb{Z}_d$.  (Here, and elsewhere, we consider integers as elements of $\mathbb{Z}_d$ via the canonical homomorphism $\mathbb{Z} \rightarrow \mathbb{Z}_d$.)

First note that, for any integer $i$ with $- l \cdot c \leq i \leq k \cdot c$, our assumption about $c$ implies $$ 1 \leq r \leq l \cdot c + i +r \leq (k+l) \cdot c +r \leq (k+l) \cdot c + \delta(d) \leq d-1,$$ and therefore, considering $$B=\{ l \cdot c + i +r \mbox{   } | \mbox{   } - l \cdot c \leq i \leq k \cdot c \}$$ as a subset of $\mathbb{Z}_d,$ we have $0 \not \in B.$

Furthermore, in $\mathbb{Z}_d$ we have $$(k-l) \cdot a = (k-l) \cdot u \cdot q =  \delta(d) \cdot q - d \cdot v \cdot q = \delta(d) \cdot q = l \cdot c +r,$$
and therefore $$kA - lA = \{(k-l) \cdot a +i \mbox{   } | \mbox{   } - l \cdot c \leq i \leq k \cdot c \}=B.$$  Since $0 \not \in B,$ $A$ is indeed $(k,l)$-sum-free in $\mathbb{Z}_d$.

Furthermore, since $c < d$, we see that $|A|=c+1$, as claimed.  
 $\quad \Box$

The lower bound in Theorem \ref{boundgeneral} now follows from Proposition \ref{delta} and the following lemma.

\begin{lem} \label{geq}
Suppose that $d$ is a divisor of $v$.  Then $$\lambda_{k,l}(G) \geq \lambda_{k,l}(\mathbb{Z}_d)  \cdot \frac{n}{d}.$$
\end{lem}

\emph{Proof.}  Since $d$ is a divisor of $v$, there is a subgroup $H$ of $G$ of index $d$ for which $$G/H \cong \mathbb{Z}_d.$$  Let $\Phi : G \rightarrow G/H$ be the canonical homomorphism from $G$ to $G/H$, and let $\Psi : G/H \rightarrow \mathbb{Z}_d$ be the isomorphism from $G/H$ to $\mathbb{Z}_d.$  Then, for any $(k,l)$-sum-free set $A \subseteq \mathbb{Z}_d$, the set $\Phi^{-1}(\Psi^{-1}(A))$ is a $(k,l)$-sum-free set in $G$ and has size $\frac{n}{d} \cdot |A|$. $\quad \Box$

\section{$(k,l)$-sum-free sets in cyclic groups}

In this section we analyze $(k,l)$-sum-free arithmetic progressions in $\mathbb{Z}_n$ and prove Theorems \ref{cyclic} and \ref{(k,l)general}.  This was carried out by Hamidoune and Plagne in \cite{HamPla:2003a} with the assumption that $n$ and $k-l$ are relatively prime; here we drop that assumption but follow their approach. 

A subset $A$ of $\mathbb{Z}_n$ is an arithmetic progression of difference $d \in \mathbb{Z}_n$, if $$A=\{a, a+d, \dots, a+ c \cdot d \}$$ for some $a \in \mathbb{Z}_n$ and non-negative integer $c$.  We let $A_{k,l}(n)$ be the set of $(k,l)$-sum-free arithmetic progression in $\mathbb{Z}_n$. 
We also let $B_{k,l}(n)$ and $C_{k,l}(n)$ be the sets of those sequences in $A_{k,l}(n)$ whose difference is not relatively prime to $n$, and relatively prime to $n$, respectively.  Note that a sequence can belong to both $B_{k,l}(n)$ and $C_{k,l}(n)$ only if it contains exactly 1 term, and that sequences in $B_{k,l}(n)$ are each contained in a proper coset in $\mathbb{Z}_n$, while no sequence in $C_{k,l}(n)$ with more than one term is contained in a proper coset.

We introduce the following notations. 
$$\alpha_{k,l}(\mathbb{Z}_n)=\max\{|A| \mbox{   } | \mbox{   } A \in A_{k,l}(n) \}$$ 
$$\beta_{k,l}(\mathbb{Z}_n)=\max\{|A| \mbox{   } | \mbox{   } A \in B_{k,l}(n) \}$$
$$\gamma_{k,l}(\mathbb{Z}_n)=\max\{|A| \mbox{   } | \mbox{   } A \in C_{k,l}(n) \}$$ 
Clearly, $\alpha_{k,l}(\mathbb{Z}_n)=\max\{\beta_{k,l}(\mathbb{Z}_n), \gamma_{k,l}(\mathbb{Z}_n) \}$.

We also let $D(n)$ be the set of all divisors of $n$ which are greater than 1.  Furthermore, we separate the elements of $D(n)$ into subsets $D_1(n)$ and $D_2(n)$ according to whether they do not or do divide $k-l$, respectively.  Then the following are clear: 

\begin{itemize}

\item $D_1(n)=\emptyset$ if, and only if, $k-l$ is divisible by $n$; 

\item $D_2(n)=\emptyset$ if, and only if, $k-l$ and $n$ are relatively prime; and

\item $D_1(n) \not =\emptyset$ and $D_2(n) \not =\emptyset$ if, and only if, $1 < \mathrm{gcd}(n,k-l) < n$.

\end{itemize}

The next three propositions summarize our results on $\alpha_{k,l}(\mathbb{Z}_n)$, $\beta_{k,l}(\mathbb{Z}_n)$, and $\gamma_{k,l}(\mathbb{Z}_n)$.  We start with $\beta_{k,l}(\mathbb{Z}_n)$.

\begin{prop} \label{beta}
The maximum size $\beta_{k,l}(\mathbb{Z}_n)$ of a $(k,l)$-sum-free arithmetic progression in $\mathbb{Z}_n$ whose difference is not relatively prime to $n$ satisfies the following.

(i) If $k-l$ is divisible by $n$, then $\beta_{k,l} (\mathbb{Z}_n)=0$.

(ii) If $k-l$ and $n$ are relatively prime, then $\beta_{k,l} (\mathbb{Z}_n)=\frac{n}{p}$ where $p$ is the smallest prime divisor of $n$.

(iii) If $1 < \mathrm{gcd}(n,k-l) < n$, then we have $$\frac{n}{\rho_1} \leq \beta_{k,l} (\mathbb{Z}_n) \leq \max \left\{ \frac{n}{\rho_1}, \frac{n}{2\rho_2} \right\},$$ where $\rho_1$
and $\rho_2$ are the smallest elements of $D_1(n)$ and $D_2(n)$, respectively.

\end{prop}

\emph{Proof.} If $n$ divides $k-l$, then for any $a \in \mathbb{Z}_n$ we have $ka=la$.  This implies (i).  Statements (ii) and (iii) will follow from the following three claims.

\emph{Claim 1.}  Suppose that $d \in D_1(n)$.  Then the set $$A=\left\{1+i \cdot d \mbox{   } | \mbox{   } 0 \leq i \leq \frac{n}{d}-1 \right\}$$ is an arithmetic progression in $B_{k,l}(n)$, has size $|A|=\frac{n}{d}$, and is $(k,l)$-sum-free.

\emph{Proof of Claim 1.}  Clearly, $A$ belongs to $B_{k,l}(n)$ and has size $|A|=\frac{n}{d}$.    Furthermore, $$kA-lA= \left\{ (k-l) + d \cdot j \mbox{   } | \mbox{   } -l \cdot \left(\frac{n}{d}-1 \right)  \leq j \leq k \cdot \left(\frac{n}{d}-1 \right)  \right\}.$$  Since $d |n$ but $d \not | (k-l)$, we have $0 \not \in kA-lA$ which means that $A$ is $(k,l)$-sum-free. 

\emph{Claim 2.}  Suppose that $H$ is a subgroup of $\mathbb{Z}_n$ of index $d$, and that $A$ is a $(k,l)$-sum-free subset of $\mathbb{Z}_n$ (not necessarily an arithmetic progression) which lies in a single coset of $H$.  Then $|A| \leq \frac{n}{d}$.

\emph{Proof of Claim 2.}  Clearly, $A \subseteq a+H$ implies $|A| \leq |H|= \frac{n}{d}$. 

\emph{Claim 3.}  Suppose again that $H$ is a subgroup of $\mathbb{Z}_n$ of index $d$, and that $A$ is a $(k,l)$-sum-free subset of $\mathbb{Z}_n$ which lies in a single coset of $H$.  If $d \in D_2(n)$, then $|A| \leq \frac{n}{2d}$.

\emph{Proof of Claim 3.}  Note that $H$ is a cyclic group of order $n/d$ and $$H=\left\{0,d,2d, \dots, \frac{n}{d}-1 \right\}.$$  Since $A$ lies in a single coset of $H$, so do $kA$ and $lA$.  But $k-l$ is divisible by $d$, so $ka-la \in H$, and therefore the sets $kA$ and $lA$ lie in the same coset of $H$.  Thus we have $$|kA \cup lA| \leq |H|=\frac{n}{d}.$$  But $A$ is $(k,l)$-sum-free, so $kA$ and $lA$ must be disjoint, hence $$|kA| + |lA| \leq \frac{n}{d}.$$  Now clearly $(k-1)a+A \subseteq kA$, so $|A| \leq |kA|$; similarly, $|A| \leq |lA|$.  This implies that $$|A| + |A| \leq \frac{n}{d}.$$  $\quad \Box$

Next, we turn to $\gamma_{k,l}(\mathbb{Z}_n)$.

\begin{prop}  \label{gamma} 
The maximum size $\gamma_{k,l}(\mathbb{Z}_n)$ of a $(k,l)$-sum-free arithmetic progression in $\mathbb{Z}_n$ whose difference is relatively prime to $n$ satisfies
$$\left \lfloor \frac{n-1-\delta}{k+l} \right \rfloor +1 \leq \gamma_{k,l}(\mathbb{Z}_n) \leq \left \lfloor \frac{n-2}{k+l} \right \rfloor +1,$$ where $\delta= \mathrm{gcd}(n,k-l)$.

\end{prop}

\emph{Proof.} The lower bound follows directly from Proposition \ref{delta}.  

For the upper bound, suppose that $d \in \mathbb{Z}_n$ and $\mathrm{gcd}(d,n)=1$, and let $a \in \mathbb{Z}_n$.  We need to show that, if the set $$A=\{a,a+d, \dots,a+c\cdot d\}$$ is $(k,l)$-sum-free in $\mathbb{Z}_n$, then $$(k+l) \cdot c \leq n-2.$$

Suppose, indirectly, that $$(k+l) \cdot c \geq n-1;$$  then we have $$\{(k-l) \cdot a +i \cdot d \mbox{   } | \mbox{   } -l \cdot c \leq i \leq k \cdot c \}   \supseteq  \{(k-l) \cdot a +j \cdot d \mbox{   } | \mbox{   } 0 \leq j \leq n-1 \}.$$  Now the left-hand side equals $kA-lA$.  Since $\mathrm{gcd}(d,n)=1$, the right-hand side equals the entire group $\mathbb{Z}_n$.  But then $kA-lA$ must contain 0, which is a contradiction. 
$\quad \Box$ 

We can now combine Propositions \ref{beta} and \ref{gamma} to get results for the maximum size of $(k,l)$-sum-free arithmetic progressions in $\mathbb{Z}_n$.

\begin{prop}  \label{alpha} 

The maximum size $\alpha_{k,l}(\mathbb{Z}_n)$ of a $(k,l)$-sum-free arithmetic progression in $\mathbb{Z}_n$ satisfies the following.

(i) If $k-l$ is divisible by $n$, then $\alpha_{k,l}(\mathbb{Z}_n)=0$.

(ii) If $k-l$ and $n$ are relatively prime, then $$\alpha_{k,l}(\mathbb{Z}_n) = \max \left \{ \frac{n}{p} ,\left \lfloor \frac{n-2}{k+l} \right \rfloor +1 \right \}$$ where $p$ is the smallest prime divisor of $n$.

(iii) If $1 < \mathrm{gcd}(n,k-l) < n$, then we have  
$$\max \left \{ \frac{n}{\rho_1} ,\left \lfloor \frac{n-1-\delta}{k+l} \right \rfloor +1 \right \} \leq \alpha_{k,l}(\mathbb{Z}_n) \leq \max \left \{ \frac{n}{\rho_1} , \frac{n}{2\rho_2} ,\left \lfloor \frac{n-2}{k+l} \right \rfloor +1 \right \},$$ where $\delta= \mathrm{gcd}(n,k-l)$, and $\rho_1$
and $\rho_2$ are the smallest elements of $D_1(n)$ and $D_2(n)$, respectively.

\end{prop}

It is easy to see that the bounds in Proposition \ref{alpha} are tight.

Now we are ready to prove Theorem \ref{cyclic}.  Due to the following result in \cite{HamPla:2003a}, our task is not difficult.

\begin{thm} [Hamidoune and Plagne, \cite{HamPla:2003a}] \label{general}

Let $\epsilon$ be 0 if $n$ is even and 1 if $n$ is odd.  Then we have the following bounds.

$$\max_{d| v} \left\{ \alpha_{k,l}(\mathbb{Z}_d) \cdot \frac{n}{d} \right\} \leq \lambda_{k,l}(G) \leq 
\max \left\{ \frac{n-\epsilon}{k+l} , \max_{d| v} \left\{ \alpha_{k,l}(\mathbb{Z}_d) \cdot \frac{n}{d} \right\} \right\}$$

\end{thm}

\emph{Proof of Theorem \ref{cyclic}.}  If $k-l$ is divisible by $n$, Theorem \ref{cyclic} obviously holds as both sides equal zero, so let's assume otherwise.  By Theorem \ref{general}, it suffices to prove that $$ \left \lfloor \frac{n-\epsilon}{k+l} \right \rfloor \leq \max_{d| n} \left\{ \alpha_{k,l}(\mathbb{Z}_d) \cdot \frac{n}{d} \right\}.$$  
By Proposition \ref{alpha}, this statement follows once we prove 
\begin{eqnarray} \label{epsilon}
 \left \lfloor \frac{n-\epsilon}{k+l} \right \rfloor \leq \max_{d| n} \left\{ \max \left \{ \frac{d}{\rho_1(d)} ,\left \lfloor \frac{d-1-\delta(d)}{k+l} \right \rfloor +1 \right \}  \cdot \frac{n}{d} \right\},
\end{eqnarray}
where $\rho_1(d)$ is the smallest divisor of $d$ which does not divide $k-l$.
 (Note that in the case when $\delta=1$, $\rho_1(d)$ is simply the smallest prime dividing $d$, thus we do not need to consider cases (ii) and (iii) of Proposition \ref{alpha} separately.)

Now $\rho_1=\rho_1(n)$ does not divide $k-l$, so we must have $\delta(\rho_1)=\mathrm{gcd}(\rho_1, k-l) < \rho_1$.  Therefore, since $\rho_1$ divides $n$, we have 

$$\max_{d| n} \left\{ \left( \left \lfloor \frac{d-1-\delta(d)}{k+l} \right \rfloor +1 \right) \cdot \frac{n}{d}  \right\} \geq 
\left( \left \lfloor \frac{\rho_1-1-\delta(\rho_1)}{k+l} \right \rfloor +1 \right) \cdot \frac{n}{\rho_1} \geq \frac{n}{\rho_1}.$$

We then have  

${\displaystyle
\max_{d| n} \left\{ \max \left \{ \frac{d}{\rho_1(d)} ,\left \lfloor \frac{d-1-\delta(d)}{k+l} \right \rfloor +1 \right \}  \cdot \frac{n}{d} \right\}  = }$  

\hspace*{1in} ${\displaystyle=  \max \left \{  \max_{d| n} \left\{ \frac{n}{\rho_1(d)} \right \}, \max_{d| n} \left\{ \left( \left \lfloor \frac{d-1-\delta(d)}{k+l} \right \rfloor +1 \right) \cdot \frac{n}{d}  \right\} \right\}} $ 

\hspace*{1in} ${\displaystyle =   \max \left \{   \frac{n}{\rho_1} , \max_{d| n} \left\{ \left( \left \lfloor \frac{d-1-\delta(d)}{k+l} \right \rfloor +1 \right) \cdot \frac{n}{d}  \right\} \right\} }$

\hspace*{1in} ${\displaystyle =   \max_{d| n} \left\{ \left(\left \lfloor \frac{d-1-\delta(d)}{k+l} \right \rfloor +1 \right) \cdot \frac{n}{d} \right \}  . }$

Therefore, (\ref{epsilon}) is equivalent to 
$$ \left \lfloor \frac{n-\epsilon}{k+l} \right \rfloor \leq \max_{d| n} \left\{ \left(\left \lfloor \frac{d-1-\delta(d)}{k+l} \right \rfloor +1 \right) \cdot \frac{n}{d} \right \}  .$$
But this inequality clearly holds, since
 \begin{eqnarray*} 
\max_{d| n} \left\{ \left(\left \lfloor \frac{d-1-\delta(d)}{k+l} \right \rfloor +1 \right) \cdot \frac{n}{d} \right \}  & \geq & \left \lfloor \frac{n-1-\delta}{k+l} \right \rfloor +1 \\
& \geq & \left \lfloor \frac{n-1-(k-l)}{k+l} \right \rfloor +1 \\ 
& = &  \left \lfloor \frac{n+(2l-1)}{k+l} \right \rfloor \\
& \geq & \left \lfloor \frac{n-\epsilon}{k+l} \right \rfloor.
\end{eqnarray*}
$\quad \Box$

\emph{Proof of Theorem \ref{(k,l)general}.}  By Theorems \ref{cyclic} and \ref{general}, here we need to show that our assumptions imply 
\begin{eqnarray} \label{epsilon1}
 \left \lfloor \frac{n-\epsilon}{k+l} \right \rfloor \leq \max_{d| v} \left\{ \max \left \{ \frac{d}{\rho_1(d)} ,\left \lfloor \frac{d-1-\delta(d)}{k+l} \right \rfloor +1 \right \}  \cdot \frac{n}{d} \right\},
\end{eqnarray}
where $\rho_1(d)$ is the smallest divisor of $d$ which does not divide $k-l$.  (The only difference between (\ref{epsilon}) and (\ref{epsilon1}) is that in (\ref{epsilon1}) only divisors of $v$ are considered.)

In a similar manner as before, we use the fact that $\rho_1(v)$ does not divide $k-l$ to conclude that the right hand side equals 
$$\max_{d| v} \left\{ \left(\left \lfloor \frac{d-1-\delta(d)}{k+l} \right \rfloor +1 \right) \cdot \frac{n}{d} \right \}  . $$

Now let $d_0$ be a divisor of $v$ which is not congruent to any integer between 1 and $\delta(d_0)$ (inclusive) (mod $k+l$).  Then the remainder of $d_0-1-\delta(d_0)$ when divided by $k+l$ is at most $k+l-1-\delta(d_0)$.  Therefore, we have
 \begin{eqnarray*} 
\max_{d| v} \left\{ \left(\left \lfloor \frac{d-1-\delta(d)}{k+l} \right \rfloor +1 \right) \cdot \frac{n}{d} \right \}  & \geq & \left(\left \lfloor \frac{d_0-1-\delta(d_0)}{k+l} \right \rfloor +1 \right) \cdot \frac{n}{d_0} \\
& \geq & \left( \frac{d_0-(k+l)}{k+l}  +1 \right) \cdot \frac{n}{d_0} \\ 
& = & \frac{n}{k+l}, 
\end{eqnarray*}
proving (\ref{epsilon1}).
$\quad \Box$

\section{$(3,1)$-sum-free sets in cyclic groups}

In this section we prove Theorem \ref{(3,1)cyclic} and find $\lambda_{3,1}(\mathbb{Z}_n)$ explicitly.  First, we evaluate $\alpha_{3,1}(\mathbb{Z}_n)$.  We note that, while Proposition \ref{alpha} (ii) readily yields
$$\alpha_{2,1}(\mathbb{Z}_n) = \left\{ 
\begin{array}{cll}
\frac{n}{2} & \mbox{   if    } & 2|n; \\ \\
\left \lfloor \frac{n+1}{3} \right \rfloor & \mbox{   if    } & 2 \not |n;
\end{array}
\right.$$
evaluating $\alpha_{3,1}(\mathbb{Z}_n)$ requires a bit more work.   

\begin{prop} \label{3,1}

The maximum size $\alpha_{3,1}(\mathbb{Z}_n)$ of a $(3,1)$-sum-free arithmetic progression in $\mathbb{Z}_n$ is given as follows: 
$$\alpha_{3,1}(\mathbb{Z}_n) = \left\{ 
\begin{array}{cll}
\frac{n}{3} & \mbox{   if    } & 3|n; \\ \\
\left \lfloor \frac{n+2}{4} \right \rfloor & \mbox{   if    } & 3 \not |n  \mbox{   and   } n  \not \equiv 2 \mbox{    \rm{(mod}    } 8);\\ \\
\frac{n-2}{4} & \mbox{   if    } & 3 \not |n \mbox{   and   } n  \equiv 2 \mbox{    \rm{(mod}    } 8).
\end{array}
\right.$$

\end{prop}

\emph{Proof.}  Let $\alpha_{3,1}(n)=\alpha$.  If $n = 2$, the claim holds, so we assume that $n \geq 3$.  We distinguish several cases.

\emph{Case 1: $2 \not | n$ and $3 \not | n$.}  In this case Proposition \ref{alpha} (ii) applies, and $$\alpha = \left \lfloor \frac{n+2}{4} \right \rfloor.$$

\emph{Case 2: $2 \not | n$ and $3  | n$.}  Proposition \ref{alpha} (ii) applies again; we get $$\alpha = \max \left \{  \frac{n}{3}, \left \lfloor \frac{n+2}{4} \right \rfloor \right\} = \frac{n}{3}.$$

\emph{Case 3: $2  | n$ and $3  | n$.}  In this case Proposition \ref{alpha} (iii) applies with $\delta=2$, $\rho_1=3$, and $\rho_2=2$; we get $$\max \left \{  \frac{n}{3}, \left \lfloor \frac{n+1}{4} \right \rfloor \right\} \leq \alpha \leq \max \left \{  \frac{n}{3}, \frac{n}{4}, \left \lfloor \frac{n+2}{4} \right \rfloor \right\},$$  which again implies $$\alpha = \frac{n}{3}.$$

\emph{Case 4: $4  | n$ and $3  \not | n$.} Again Proposition \ref{alpha} (iii) applies -- this time with $\delta=2$, $\rho_1=4$, and $\rho_2=2$.  Therefore we get $$\max \left \{  \frac{n}{4}, \left \lfloor \frac{n+1}{4} \right \rfloor \right\} \leq \alpha \leq \max \left \{  \frac{n}{4}, \left \lfloor \frac{n+2}{4} \right \rfloor \right\},$$  which gives $$\alpha = \frac{n}{4}.$$ 

\emph{Case 5: $n \equiv 2 $ (mod 4) and $3  \not | n$.} Again Proposition \ref{alpha} (iii) applies --- this time with $\delta=2$, $\rho_1 \geq 5$, and $\rho_2=2$.  Therefore we get $$\max \left \{  \frac{n}{\rho_1}, \left \lfloor \frac{n+1}{4} \right \rfloor \right\} \leq \alpha \leq \max \left \{  \frac{n}{\rho_1}, \frac{n}{4}, \left \lfloor \frac{n+2}{4} \right \rfloor \right\},$$  which yields only $$\alpha \in \left \{ \frac{n-2}{4}  ,  \frac{n+2}{4}  \right \}.$$ 

To continue further, we separate the cases of $n \equiv 2 $ (mod 8) and $n \equiv 6 $ (mod 8).  

\emph{Case 5.1.}  Let us first consider the case when $n \equiv 6 $ (mod 8).  With $a = \frac{n+2}{8}$ and $c=\frac{n-2}{4}$, we let $$A=\{a, a+1, \dots, a+c \}.$$  Then 
$$
3A-A  =  \{ 2a -c +i \mbox{   } |\mbox{   } 0 \leq i \leq 4c \} 
  =  \{ 1 + i \mbox{   } |\mbox{   } 0 \leq i \leq n-2 \}
= \mathbb{Z}_n \setminus \{0 \},$$       
so $A$ is (3,1)-sum-free in $\mathbb{Z}_n$ of size $c+1=\frac{n+2}{4}$.

\emph{Case 5.2.}  Now suppose that $n \equiv 2 $ (mod 8).  We prove that $\alpha = \frac{n-2}{4}$.  Suppose, indirectly, that $\alpha = \frac{n+2}{4}$ and there is a (3,1)-sum-free arithmetic progression $$A=\{a, a+d, \dots, a+c \cdot d \}$$ in $\mathbb{Z}_n$ of size $c+1=\frac{n+2}{4}$.  Similarly to above, $$
3A-A  =  \{ 2a -c \cdot d +i \cdot d \mbox{   } |\mbox{   } 0 \leq i \leq 4c \} 
  =  \{ 2a -c \cdot d +i \cdot d \mbox{   } |\mbox{   } 0 \leq i \leq n-2 \}.$$

By Proposition \ref{beta} (iii), we have $$\beta_{3,1}(n) \leq \max \left \{\frac{n}{\rho_1}, \frac{n}{4}  \right \} =\frac{n}{4} ; $$  so we have $\beta_{3,1}(n) < \alpha$.  Therefore, we must have $\mathrm{gcd}(d,n)=1$, which implies that $$|3A-A|=n-1.$$   
Since $A$ is (3,1)-sum-free, $0 \not \in 3A-A$, and this can only occur if $$2a -c \cdot d +(n-1) \cdot d \equiv 0 \mbox{    \rm{(mod}    } n).$$  A simple parity argument provides a contradiction: $2a -c \cdot d +(n-1) \cdot d$ is odd, so it cannot be divisible by $n$.  $\quad \Box$

\emph{Proof of Theorem \ref{(3,1)cyclic}.}  As previously, we let $D(n)$ be the set of divisors of $n$ which are greater than 1.  We introduce the following six (potentially empty) subsets of $D(n)$, as well as some notations.

\begin{tabular}{lllllll}
$E_1(n)$ & = & $\{ d \in D(n) \mbox{   } | \mbox{   } 3|d \}$ & \mbox{    } &
$ e_1 $ & = & $\max_{d \in E_1(n)} \left\{ \frac{d}{3} \cdot \frac{n}{d} \right\}$ \\ \\
$E_2(n)$ & = & $\{ d \in D(n) \mbox{   } | \mbox{   } d \equiv 3 (4), 3 \not |d \}$ &  \mbox{    } &
$e_2 $& = & $\max_{d \in E_2(n)} \left\{ \frac{d+1}{4} \cdot \frac{n}{d} \right\}$ \\ \\
$E_3(n)$ & = & $\{ d \in D(n) \mbox{   } | \mbox{   } 4|d, 3 \not |d \}$  & \mbox{    } & 
$e_3$ & = & $\max_{d \in E_3(n)} \left\{ \frac{d}{4} \cdot \frac{n}{d} \right\}$ \\ \\
$E_4(n)$ & = & $\{ d \in D(n) \mbox{   } | \mbox{   } d \equiv 1 (4), 3 \not |d \}$ &  \mbox{    } &
$e_4 $& = & $\max_{d \in E_4(n)} \left\{ \frac{d-1}{4} \cdot \frac{n}{d} \right\}$ \\ \\
$E_5(n)$ & = & $\{ d \in D(n) \mbox{   } | \mbox{   } d \equiv 6 (8), 3 \not |d \}$ &  \mbox{    } &
$e_5 $& = & $\max_{d \in E_5(n)} \left\{ \frac{d+2}{4} \cdot \frac{n}{d} \right\} $ \\ \\
$E_6(n)$ & = & $\{ d \in D(n) \mbox{   } | \mbox{   } d \equiv 2 (8), 3 \not |d \}$ &  \mbox{    } &
$e_6$ & = & $\max_{d \in E_6(n)} \left\{ \frac{d-2}{4} \cdot \frac{n}{d} \right\}$
\end{tabular}

\noindent (We have the understanding that $\max \emptyset = 0.$)

Then we have $$D(n)= \cup_{i=1}^6 E_i(n);$$ furthermore, by Theorem \ref{cyclic} and Proposition \ref{3,1}, we have $$\lambda_{3,1}(\mathbb{Z}_n) = \max \{e_i | 1 \leq i \leq 6 \}.$$

For any $i \in \{1,2, \dots, 6\}$ for which $E_i(n) \not = \emptyset$, we let $$p_i= \min \{E_i(n) \}$$ and $$n_i= \max \{E_i(n) \}.$$ 

Now suppose that $E_5(n) \not = \emptyset$.  Then $E_2(n) \not = \emptyset$, and $p_5=2 \cdot p_2$.  Therefore 
$$e_5=\frac{p_5+2}{4} \cdot \frac{n}{p_5}=\frac{p_2+1}{4} \cdot \frac{n}{p_2}=e_2.$$
We can similarly show that, if $E_6(n) \not = \emptyset$, then $E_4(n) \not = \emptyset$ and $e_6=e_4.$  Therefore, we see that $$\lambda_{3,1}(\mathbb{Z}_n) = \max \{e_i | 1 \leq i \leq 4 \}.$$

Next, observe that, if $E_i(n) \not = \emptyset$ for some $i \in \{1,2, 3\}$, then $e_i \geq e_{j}$ for all $i < j \leq 4$.

Now we consider the following cases.

\emph{Case 1.}  Suppose that $n$ has divisors which are congruent to 3 mod 4, and let $p$ be the smallest such divisor.  If $p=3$, then $E_1(n) \not = \emptyset$, thus $$\lambda_{3,1}(\mathbb{Z}_n) = e_1 = \frac{n}{3}.$$  If, on the other hand, $p>3$, then $E_1(n) = \emptyset$ but $E_2(n) \not = \emptyset$, thus $$\lambda_{3,1}(\mathbb{Z}_n) = e_2 = \frac{p+1}{p} \cdot \frac{n}{4}.$$  

\emph{Case 2.}  Suppose that $n$ has no divisors which are congruent to 3 mod 4, but that $n$ is divisible by 4.  In this case, $E_1(n) = E_2(n) = \emptyset$ but $E_3(n) \not = \emptyset$, thus $$\lambda_{3,1}(\mathbb{Z}_n) = e_3 = \frac{n}{4}.$$ 

\emph{Case 3.}  Suppose that $n$ has no divisors which are congruent to 3 mod 4, and that $n$ is not divisible by 4.  In this case, $E_1(n) = E_2(n) = E_3(n)= \emptyset$ but $E_4(n) \not = \emptyset$, thus $$\lambda_{3,1}(\mathbb{Z}_n) = e_4 = \frac{n_4 -1}{4} \cdot \frac{n}{n_4}.$$  If $n$ is odd, then $n_4=n$; if $n$ is even, then (since $n$ is not divisible by 4), $n_4=\frac{n}{2}$.  In either case, we get 
$$\lambda_{3,1}(\mathbb{Z}_n) = e_4 = \frac{n_4 -1}{4} \cdot \frac{n}{n_4}= \left \lfloor \frac{n}{4} \right \rfloor.$$

The claims of Theorem \ref{(3,1)cyclic} now readily follow.
$\quad \Box$

\section{Further comments and open questions}

In this final section, we discuss some interesting open questions.

Our first question is about a possible generalization of Theorems \ref{(2,1)cyclic} and \ref{(3,1)cyclic}.  Note that, according to Theorem \ref{boundgeneral}, we have $$\lambda_{k,l}(\mathbb{Z}_n) \leq \max_{ d |n} \left \{ \left( \left \lfloor \frac{d-2}{k+l}  \right \rfloor   +1 \right) \cdot \frac{n}{d} \right \}.$$

\begin{que} \label{one} Let $D(n)$ be the set of divisors of $n$ (which are greater than 1).   
Given distinct positive integers $k$ and $l$, is there a subset $D_{k,l}(n)$ of $D(n)$ so that $$\lambda_{k,l}(\mathbb{Z}_n) = \max_{ d \in D_{k,l}(n)} \left \{ \left( \left \lfloor \frac{d-2}{k+l}  \right \rfloor   +1 \right) \cdot \frac{n}{d} \right \}?$$ 
\end{que}
  
As we see from (\ref{relprime}), Question \ref{one} holds with $D_{k,l}(n)=D(n)$ when $n$ and $k-l$ are relatively prime, in particular, for sum-free sets.  According to Theorem \ref{(3,1)cyclic}, the set $$D_{3,1}(n) = \{ d \in D(n) | d \not \equiv 2 \mbox{   } (\mathrm{  mod  } \mbox{   }  4) \}$$ works for $(k,l)=(3,1)$.  (Note that, if it exists, $D_{k,l}(n)$ is not necessarily unique.)

Moving on to general abelian groups, we observe that, by Lemma \ref{geq}, we have $$\lambda_{k,l}(G) \geq \lambda_{k,l}(\mathbb{Z}_{v}) \cdot  \frac{n}{v}.$$  Then one of course wonders the following.

\begin{que} \label{two}  Given distinct positive integers $k$ and $l$, is $$\lambda_{k,l}(G)= \lambda_{k,l}(\mathbb{Z}_{v}) \cdot  \frac{n}{v}?$$
\end{que}

According to Theorem \ref{cyclic}, Question \ref{two} is equivalent to asking:  is 
$$\lambda_{k,l}(G) = \max_{d| v} \left\{ \alpha_{k,l}(\mathbb{Z}_d) \cdot \frac{n}{d} \right\}?$$

Note that Theorem \ref{(2,1)all} of Green and Ruzsa affirms Question \ref{two} for sum-free sets.  Theorem \ref{(k,l)general} exhibits some other cases when the equality also holds.  In particular, as a consequence of Theorem \ref{(k,l)general}, we see that  
$$\lambda_{3,1}(G) = \lambda_{3,1}(\mathbb{Z}_{v}) \cdot  \frac{n}{v}$$ holds when $v$ (iff $n$) has at least one prime divisor which is congruent to 3 mod 4, or when $v$ is divisible by 4.  So the only cases left open are when $v=P$ or $v=2P$ where $P$ is the product of primes all of whom are congruent to 1 mod 4.

Next, we are interested in characterizing all $(k,l)$-sum-free subsets of maximum size.

\begin{que} \label{three}  What are the $(k,l)$-sum-free subsets $A$ of $G$ with size $|A|=\lambda_{k,l}(G)?$
\end{que}

A pleasing answer is given by Bier and Chin \cite{BieChi:2001a} for the case when $k \geq 3$ and $G \cong \mathbb{Z}_p$ where $p$ is an odd prime: in this case $A$ is an arithmetic progression.  The same answer was given by Diananda and Yap \cite{DiaYap:1969a} earlier for the case when $(k,l)=(2,1)$ (that is, when $A$ is sum-free) and $G \cong \mathbb{Z}_p$ with $p$ not congruent to 1 mod 3; however, for $p=3m+1$ the set $$A=\{m, m+2, m+3, \dots, 2m-1, 2m+1 \}$$ is also sum-free with maximum size.  More generally, the answer to Question \ref{three} is known for $(k,l)=(2,1)$ and when $n$ has at least one divisor not congruent to 1 mod 3: in this case $A$ is the union of arithmetic progressions of the same length.  More precisely, there is a subgroup $H$ in $G$ so that $G/H$ is cyclic and $$A=\{(a+H) \cup (a+d+H) \cup \cdots \cup (a+c \cdot d+H)\}$$ for some $a, d \in G$ and integer $c$.  These and other results can be found in \cite{WalStrWal:1972a}.

More ambitiously, one may ask for a characterization of all ``large'' (but not necessarily maximal) $(k,l)$-sum-free sets in $G$.  Can one, for example, describe explicitly all $(k,l)$-sum-free sets of size greater than $n/(k+l)$?  Hamidoune and Plagne \cite{HamPla:2003a} carry this out for sum-free sets of size at least $n/3$ in arbitrary groups.  Other results can be found in the papers of Davydov and Tombak \cite{DavTom:1989a} and Lev \cite{Lev:2005a}, \cite{Lev:2006a}.    

Our final question is about the number of $(k,l)$-sum-free subsets in $G$, which we here denote by $N_{k,l}(G)$.

\begin{que} \label{four}  What is the cardinality $N_{k,l}(G)$ of the set of $(k,l)$-sum-free subsets in $G$?
\end{que}

Clearly, any subset of a $(k,l)$-sum-free set is also $(k,l)$-sum-free, so the answer to Question \ref{four} is at least $$N_{k,l}(G) \geq 2^{\lambda_{k,l}(G)}.$$  But there are indications that the number is not much larger.  In fact, for sum-free sets we have the following result of Green and Ruzsa \cite{GreRuz:2005a}:
$$N_{2,1}(G) = 2^{\lambda_{2,1}(G)+o(1)n},$$ where $o(1)$ approaches zero as $n$ goes to infinity.  They have a more accurate approximation for the case when $n$ has a prime divisor which is congruent to 2 mod 3.  (This result had been established for even $n$ earlier by Lev, {\L}uczak, and Schoen \cite{LevLucSch:2001a} and independently by Sapozhenko \cite{Sap:2002a}.)

In closing, we mention that the analogues of our questions about the maximum size, the structure, and the number of $(k,l)$-sum-free sets (especially sum-free sets) have been investigated in non-abelian groups (see Kedlaya's papers \cite{Ked:1997a} and \cite{Ked:1998a}) and, more extensively, among the positive integers (see the works of Alon \cite{Alo:1991a}, Bilu \cite{Bil:1998a}, Calkin \cite{Cal:1990a}, Calkin and Taylor \cite{CalTay:1996a}, Cameron \cite{Cam:1987a}, Cameron and Erd\H{o}s \cite{CamErd:1990a} and \cite{CamErd:1999a}, and {\L}uczak and Schoen \cite{LucSch:2000a}).  General background references on related questions include Nathanson's book \cite{Nat:1996a}, Guy's book \cite{Guy:2004a}, and Ruzsa's papers \cite{Ruz:1993a} and \cite{Ruz:1995a}; see also \cite{Baj:2004a} and \cite{BajRuz:2003a}.

\end{document}